# Adaptive complexity regularization for linear inverse problems

**Jean-Michel Loubes**
**Carenne Ludeña**

**Abstract:** We tackle the problem of building adaptive estimation procedures for ill-posed inverse problems. For general regularization methods depending on tuning parameters, we construct a penalized method that selects the optimal smoothing sequence without prior knowledge of the regularity of the function to be estimated. We provide for such estimators oracle inequalities and optimal rates of convergence. This penalized approach is applied to Tikhonov regularization and to regularization by projection.

**Introduction**

We are interested in recovering an unobservable signal $x_0$ based on observations

$$y(t_i) = Tx_0(t_i) + \varepsilon_i, \tag{1}$$

where $T : X \to Y$ is a known compact linear operator between Hilbert spaces $X, Y$. We cannot observe our target function $x_0$ directly, but rather through a blurred (by the linear filter $T$), noise corrupted sample of $y$ over a collection of discrete observation points $t_i, i = 1, \ldots, n$. Throughout the paper, we shall denote $\mathbf{y} = (y(t_i))_{i=1}^n$.

We assume that the observations $y(t_i) \in \mathbb{R}$ and that the observation noise $\varepsilon_i$ are i.i.d. realizations of a certain random variable $\varepsilon$.

We will say the problem is ill-posed if the generalized inverse operator $T^+ : Y \to X$ defined by $(T^*T)^{-1}T$, where $T^*$ stands for the transpose operator, is unbounded. We In this case, this inverse operator needs to be, in some sense, regularized.

Regularization methods replace an ill-posed problem by a family of well-posed problems. Their solution, called regularized solutions, are used as approximations of the desired solution of the inverse problem. These methods always involve some parameter measuring the closeness of the regularized and the original (unregularized) inverse problem. Rules (and algorithms) for the choice of these regularization parameters as well as convergence properties of the regularized solutions are central points in the theory of these methods.

The statistical problem has been extensively studied, although in general efficient regularization-parameter choice is still under active research. Two main kinds of estimators have been considered in the literature. First regularized estimators such as Tikhonov type estimators, and second non linear thresholded estimators. The first approach has been studied in great detail. An interesting early survey of this topics are provided in [20, 19]. In this setting, the main issues are, what kind of regularization functional should be considered, and, closely related, what the relative weight of the selected regularization functional should be. More recently, Mair and Ruymgaart in [18] or Bissantz et al. in [3] studied different regularized inverse problems and proved the optimality of the rate of convergence for their estimators assuming the regularity of the target function $x_0$ to be known. Special attention has been devoted in this setting when considering a Singular value decomposition (SVD) of operator $T$. We cite the recent work in this direction in [7, 6]. The second approach has its most popular version in the *wavelet-vaguelet* decomposition introduced in [9]. In this case the main issue is finding an appropriate basis over which $T^+$, the generalized inverse, is almost diagonal. This idea is further developed in [12] who introduce *mirror wavelets*. Closely related, Cohen et al. in [8] construct an adaptive thresholded estimator based on Galerkin's method. We also cite the recent work in [14], where they combine an SVD approach with a thresholding technique over a certain new basis.

Reconstructing the unknown target function $x_0$ is related to four issues. To the filter $T$, to the probabilistic structure of the noise, to the fact we are only observing $T(x_0)$ over a finite observation scheme and finally to the regularity of $x_0$. In order to unify notation, our assumptions will be presented in terms of an underlying basis of $Y$, $\{\psi_j\}_{j \in \mathbb{N}}$, and the increasing sequence of approximating spaces $Y_m := <\psi_j>_{j \in d_m}$. The ill-posedness of $T$ will then be defined in terms of these subspaces $Y_m$ and the discrete observation scheme. As is usual in the numerical literature, the regularity of $x_0$ will be defined in terms of operator $T$ (for a detailed discussion see [11] or [13]).

Our goal in this article is to develop algorithms providing estimators of $x_0$ that achieve optimal rates of convergence within the regularization method when the smoothness of the true solution is not known a priori. For this, we focus on model selection techniques for regularization procedures. As in the work of [5], [15] or [16], we choose the best regularization scheme among a set of regularization operators by minimizing a contrast and a well chosen penalty. Since the choice of a penalty is crucial, we provide a very general way of calibrating the penalty with respect to the regularization operators. This enables us to build optimal estimators, within the chosen regularization method, for inverse problems for two general classes of estimators, Tikhonov estimators and projection estimators.

The article is divided into four main parts. In Section 1 we describe the general framework and provide several assumptions. In Section 2 we build a general penalty which leads to an oracle inequality for regularization operators and shows optimality of the adaptive procedures. In Section 3 we apply the above results to Tikhonov regularization operators and projection operators. Section 4 is devoted to technical lemmas providing deviation bounds which are used in order to obtain non asymptotic bounds for the quadratic risk of our estimator.

## 1. Model description and general assumptions

In this section we introduce our general notation and assumptions.

We want to estimate a function $x_0 : \mathbb{R} \to \mathbb{R}$ based on the blurred and noisy observations

$$y(t_i) = Tx_0(t_i) + \varepsilon_i, \; : i = 1, \ldots, n.$$

It is important to stress that the observations depend on a fixed design $(t_1, \ldots, t_n) \in \mathbb{R}^n$. This will require introducing an empirical norm based on this design. Let $Q_n$ be the empirical measure of the co-variables $Q_n = \frac{1}{n}\sum_{i=1}^{n} \delta_{t_i}$, where $\delta$ is the Dirac function. The $L_2(Q_n)$-norm of a function $y \in Y$ is then given by

$$\|y\|_n = \left(\int y^2 dQ_n\right)^{1/2},$$

and the empirical scalar product by $<y, \varepsilon>_n = \frac{1}{n}\sum_{i=1}^{n} \varepsilon_i y(t_i)$. Remark that this empirical norm is defined over the observation space $Y$. Over the solution space $X$ we will consider the norm given by the Hilbert space structure. For sake of simplicity, we will write $\|.\|_X = \|.\|$ when no confusion is possible. Over a finite dimensional space, the norm $\|.\|$ will always stand for the Euclidean norm and if $v \in \mathbb{R}^d$, $v^t$ will stand for the transpose vector. Likewise for any matrix $A \in \mathbb{R}^{d\times r}$, $A^t$ will stand for the transpose matrix and $A^+ := (A^t A)^{-1} A^t$ for the generalized inverse. Considered as an operator, we will write $A^*$ for the adjoint of the corresponding operator.

We also introduce certain standard assumptions on the observation noise

**AN moment condition for the errors**

$\varepsilon$ is a centered random variable satisfying the moment condition $\mathbb{E}(|\varepsilon|^q/\sigma^q) \leq q!/2$ for all $q \geq 1$, with $\mathbb{E}(\varepsilon^2) = \sigma^2$.

As usual in statistics, assume that $X$ satisfies a certain smoothness condition. In this paper, we assume the following source assumption encountered typically in the inverse problems literature.

**SC source condition**

There exists $\nu > 0$ such that $x_0 \in Range((T^*T)^\nu) := \mathcal{R}((T^*T)^\nu)$

Moreover consider

$$A_{\nu,\rho} = \{x \in X, x = (T^*T)^\nu \omega, \|\omega\| \leqslant \rho\}$$

where $0 \leqslant \nu \leqslant \nu_0$, $\nu_0 > 0$ and use the further notation

$$A_\nu = \bigcup_{\rho > 0} A_{\nu,\rho} = \mathcal{R}((T^*T)^\nu) \qquad (2)$$

These sets are usually called source sets, $x \in A_\nu$ is said to have a source representation.

**Remark 1.1.** Such condition is usual in analysis of inverse problems (see [11]). It links the decay of the eigenvalues of the operator with the decay of the coefficients of the decomposition of the function in the SVD basis. Thus, the parameter $\nu$ can be seen as a smoothness parameter providing a restriction over the regularity of the function to be recovered. Indeed, let $(\lambda_j, \phi_j, \varphi_j)_{j \geqslant 1}$ be the singular value decomposition of the operator $T$, in the sense that the following decomposition holds

$$T^*Tx = \sum_{j=1}^{\infty} \lambda_j^2 <x, \phi_j> \phi_j.$$

Hence, $x_0 \in A_\nu$ if and only if

$$\sum_{j \geqslant 1} \frac{|<x_0, \phi_j>|^2}{\lambda_j^{2+4\nu}} < +\infty.$$

The parameter $\nu$ is not known a priori, hence adaptation results will be provided with respect to this smoothness parameter.

Estimating over all $X$ is in general not possible because we can only observe $Tx_0$ over the fixed design $(t_1, \ldots, t_n)$. Thus we assume that we are equipped with a sequence of nested linear subspaces whose union is dense in $Y$, $Y_1 \subset Y_2 \ldots \subset Y_m \ldots \subset Y$. We assume $dim(Y_m) = d_m$. We are interested in a subcollection of these spaces generated by a set of indices $\mathcal{M}_n$. In this paper, we will use these approximation spaces as projection spaces in order to study the data. So, denote the projection of any space $W$ over any subspace $Z$ by $\Pi_Z W$. Let $\Pi_{Y_m}^n$ stands for the projection in the empirical norm. Set also the corresponding projected operator $T_m = \Pi_{Y_m}^n T$.

Using a sieve of the space $Y$, we consider the corresponding approximation spaces in the space $X$, defined as $X_m = T_m^* Y$. By construction

$$\Pi_{X_m} = (\Pi_{Y_m}^n T)^+ \Pi_{Y_m}^n T.$$

We point out that both $T_m$ and its adjoint operator $T_m^*$ depend on the observation sequence $t_i$. However, we will usually drop this fact from the notation. To illustrate this assertion, consider the following example.

**Example 1.1.** *if $Y_m$ is generated by some orthonormal basis $\phi = (\phi_1, \ldots, \phi_{d_m})$, with respect to the $L^2$ norm over $Y$, and $T = \text{Id}$, then*

$$\Pi_{Y_m}^n y = \sum_{j=1}^{d_m} y_{j,n} \phi_j,$$

*where $y_{j,n} = <\Pi_{Y_m}^n y, \phi_j>_n$ are the solution to the projection problem under the empirical measure $Q_n$. Set $G_m = (\phi_j(t_i))_{i,j}$, $i = 1, \ldots, n$ and $j = 1, \ldots, d_m$. Thus, we may write in matrix notation*

$$\Pi_{Y_m}^n \mathbf{y} = (G_m^t G_m)^{-1} G_m^t (y(t_1), \ldots, y(t_n)).$$

An important issue is that, as above, we can always define $T_m$ in matrix notation and thus $T_m^t \mathbf{y}$ always makes sense. Moreover, if $u \in Y$, we will use indistinctly $T_m^t \mathbf{u} \in \mathbb{R}^{d_m}$ and $T_m^* u \in X_m$, the latter in operator notation.

Now, define $p$ the degree of ill-posedness of the forward operator $T$. In our case we will relate this parameter to the approximation properties of the spaces $Y_m$, i.e with the projection operator $\Pi_{Y_m}^n$ as follows.

**IP ill-posedness of the operator** Let $\mathcal{M}_n$ be an index set. For $m \in \mathcal{M}_n$ there exists a parameter $p > 0$ such that

$$\gamma^{(m)} := \|(I - \Pi_{Y_m}^n)T\| = O(d_m^{-p}).$$

$p$ is denoted the index of ill-posedness of the forward operator. In the following, this parameter is supposed to be known.

To illustrate this condition we include the following example

**Example 1.2.** *The above assumption can be seen to hold under certain conditions over operator $T$ and matrix $G_m$ defined in example (1.1). Let $(\lambda_j, \phi_j, \varphi_j)_j$ be the singular value decomposition of the forward operator $T$ and assume that there exists $p > 0$ such that $\lambda_j = O(j^{-p})$. This parameter is in the statistical literature often considered as an index of ill-posedness, since the difficulty of the estimation in inverse problem usually comes from the difficulty to invert the operator due to the decrease of its eigenvalues. But, if we let $Y_m$ be the linear subspace generated by $\{\varphi_j\}_{1 \leq j \leq d_m}$, then assume that the fixed observation design $t_i, i = 1, \ldots, n$ is such that this basis is also orthogonal in the empirical norm. Assume also that $\sup_{j=1,\ldots,d_m} \|\varphi_j\|_\infty < \infty$. Then,*

$$I - \Pi_{Y_m}^n = I - \Pi_{Y_m} + \Pi_{Y_m}^n(I - \Pi_{Y_m}),$$

*where*

$$\|I - \Pi_{Y_m}\| = O(d_m^{-p})$$

*and*

$$\begin{aligned}\|\Pi_{Y_m}^n(I - \Pi_{Y_m})u\| &= \|\Pi_{Y_m}^n(I - \Pi_{Y_m})u\|_n \\ &\leq \|(I - \Pi_{Y_m})u\|_n \leq \sup_{j=1,\ldots,j_m} \|\varphi_j\|_\infty \|(I - \Pi_{Y_m})u\|,\end{aligned}$$

*leading to the approximation property provided in [**IP**]. If the forward operator acts in a smoothing way as integrating $p$ times, for example, the eigenvalues $\lambda_j$ satisfy the required decay rates (see [13]).*

Define

$$\nu_m := \|T_m^+ \Pi_{Y_m}^n\|.$$

This quantity controls the amplification of the observation error over the solution space $X_m$. Consider

$$\gamma_m := \inf_{v \in Y_m, \|v\|=1} \|T_m^* v\|,$$

which expresses the effect of operator $T_m^*$ over the approximating subspace $Y_m$. We have as in [13], $\nu_m \geq \gamma_m$.

On the other hand this term is related to the goodness of the approximation scheme. Following the proof in [13], it can be seen that

$$\gamma_{m+1} \leq \|T^*(I - \Pi_{Y_m}^n)\| = \|(I - \Pi_{Y_m}^n)T\|.$$

The next assumption requires that $\gamma_m$ and $\gamma^{(m)}$ are of the same order, which will be written $\gamma_m \sim \gamma^{(m)}$.

**AS amplification error**

We assume

$$\gamma_m = O(d_m^{-p}). \tag{3}$$

Moreover assume there exists a positive constant $U$ such that

$$\frac{\gamma^{(m)}}{\gamma_m} \leq \sqrt{U}.$$

In the estimation procedure, we will first project the data onto a large enough approximation space indexed by $m_0$, to be selected later in this paper. For this fixed $m_0$, let $(\lambda_j, \phi_j, \varphi_j) j = 1, \ldots, d_{m_0}$ be the singular value decomposition of operator $T_{m_0}$. For any $u \in Y$ we can write $T_{m_0}^* u = \sum_{j=1}^{d_{m_0}} \lambda_j \phi_j <u, \varphi_j>_n$, which depends only on $\mathbf{u} = (u(x_1), \ldots, u(x_n))^t$. As previously recall that $G_{m_0} \in M_{d_{m_0}, n}$ is defined by $G = (\varphi_j(x_i))_{j,i}$, $i = 1, \ldots, n$, $j = 1, , d_{m_0}$. Thus, abusing notation we may write $T_{m_0}^t = D(G_{m_0} G_{m_0}^t)^{-1} G_{m_0} : \mathbb{R}^n \to \mathbb{R}^{d_{m_0}}$, where $D = D(\lambda_j)_{j=1,\ldots,d_{m_0}}$ is the diagonal matrix with entries $\lambda_j$. Since $T_{m_0}^t \mathbf{u} = T_{m_0}^* u$ the latter in operator notation, both interpretations will be used indistinctly. On the other hand, for $x \in X_{m_0}$, identified with a $d_{m_0}$ dimensional vector, we can think of $(T_{m_0} x(t_1), \ldots, T_{m_0} x(t_n)) = G_{m_0}^t D x$. So that in matrix notation also $T_{m_0}^* T_{m_0} = D^2$.

We also introduce the following assumptions

**SV** There exist positive constants $k_1 < k_2$ such that $k_1 j^{-p} \leq \lambda_j \leq k_2 j^{-p}$.
**SF** Let $\nu_j$, $j = 1, \ldots, d_{m_0}$ be the eigenvalues of matrix $G^t G$, then there exist constants $a_1 < a_2$ such that $a_1 n \leq \nu_j \leq a_2 n$.

**Remark 1.2.**

- Assumption **AS** thus establishes that the worst amplification of the error over $X_m$ is roughly equivalent to the best approximation over $Y_m$ in the empirical norm. This yields a uniform bound on the operator $TT_m^+ \Pi_{Y_m}^n$, so that regularization by $T$ compensates the bad condition of the approximation $T_m^+ \Pi_{Y_m}^n$ to $T^+$ (see [11] or [13] for further comments). Notice $T_m^+ \Pi_{Y_m}^n T$ is just the projection operator over $X_m$ so it is a bounded operator.
- Assumption **SV** is slightly stronger than **AS** as it establishes the exact order of the $\gamma_m$. It is seen to hold, for instance, in example (1.2). Assumption **SF** is necessary to assure convergence results further on. It holds also in example (1.2).

- We point out that the Assumptions are met as soon as the approximation spaces are constructed as in example (1.2), i.e called truncated SVD where the spaces are defined as $Y_m = \text{span}\{\varphi_j,\ j = 1, \ldots, d_m\}$, where we recall that $\varphi_j$ are eigenfunctions corresponding to the eigenvalues of $T^*T$ sorted in decreasing order. Since eigenfunctions are often hard to obtain explicitly, for practical purpose we are more interested in the usual piecewise polynomial projection bases or wavelet bases, though, as in [13].

## 2. Adaptive choice for regularized estimator

Consider a class of regularized estimators built using a projection and a regularization procedure. Namely let $Y_{m_0}$ be a big enough space in the sense that $m_0$ is such that

$$\|(I - \Pi_{X_{m_0}})x_0\| \leq \inf_{m \in \mathcal{M}_n}[\|(I - \Pi_{X_m})x_0\| + \sqrt{\frac{d_m}{n}}\frac{1}{\gamma_m}].$$

This quantity can be chosen so as not to depend on the unknown regularity of the solution $x_0$. Under assumption **SC** the above inequality is satisfied if the dimension of the set is such that

$$d_{m_0}^{2\nu p} \geq n^{\frac{2\nu p}{4\nu p + 2p + 1}}.$$

Thus it is enough to choose $m_0$ such that $d_{m_0} \geq n^{1/(2p+1)}$.

For $\mathcal{K}_n$ a set of indices, consider $\{\tilde{R}_k,\ k \in \mathcal{K}_n\}$ a collection of regularization operators which depend on different values of the smoothing parameters. For instance consider Tikhonov regularization operators which rely on the choice of a smoothing sequence, Landweber iteration operators which rely on the choice of a stopping index, or other general smoothing operators described in [10]. Consider the corresponding estimators

$$\hat{x}_k := \tilde{R}_k \Pi_{Y_{m_0}}^n y = R_k y, \tag{4}$$

where we have written $R_k := \tilde{R}_k \Pi_{Y_{m_0}}^n$. The behavior of such general estimators depends on the choice of the regularization sequence. From the theory of inverse problems, we know that it is possible to choose a regularization operator for which the corresponding estimator achieves the optimal rate of convergence, but this choice depends on $\nu$ defined in **SC**, which characterizes the regularity of the solution.

Our aim is building a method that picks, according to the data, an optimal $R_k$, among all the $R_k,\ k \in \mathcal{K}_n$ in such a way that optimal rates are maintained. This choice must also not depend on a priori regularity assumptions. We point out that selecting the optimal smoothing parameter in a collection of sequences, belongs to model selection theory since it is equivalent as selecting a good model among a collection of sets.

For this consider the following penalized procedure. For a given constant $r > 2$ and weights $L_k$, $k \in \mathcal{K}_n$ to be chosen, define the penalty as

$$\text{pen}(k) := r\sigma^2(1+L_k)[Tr(R_k^t R_k) + \rho^2(R_k)],$$

where $Tr(R_k^t R_k)$ is the trace and $\rho(R_k^t R_k) = \rho^2(R_k)$ is the spectral radius. Finally $\hat{k}$ is selected as the solution of

$$\hat{k} := \arg\min_{k \in \mathcal{K}_n} \left\{ \|R_k(y - T(\hat{x}_k))\|^2 + \text{pen}(k) \right\}, \quad (5)$$

which defines the estimator $\hat{x}_{\hat{k}} = R_{\hat{k}} y$. Let $x_k = R_k T x_0$ be the regularized true function, which measures the accuracy of the estimation procedure without observation noise. The following result states the asymptotic behaviour of the estimator $\hat{x}_{\hat{k}}$.

**Theorem 2.1.** *Assume* **AN** *and* **SF** *are satisfied. There exists a constant $C$ which depends on $r$ and on $T$, such that the following inequality holds true*

$$\mathbf{E}\|\hat{x}_{\hat{k}} - x_0\|^2 \leq 2\|(I - \Pi_{X_{m_0}})x_0\|^2 + C \inf_{k \in \mathcal{K}_n} \left[ \|x_k - x_0\|^2 + 2\text{pen}(k) \right] + \frac{\Sigma(d)}{n}, \quad (6)$$

*where we have set*

$$\Sigma(d) = \sum_{k \in \mathcal{K}_n} 2 \left[ \sqrt{\frac{dTr(R_k^t R_k)}{\rho^2(R_k)}} + 1 \right] \left[ \frac{d}{\rho^2(nR_k)} \right]^{-1} e^{-\sqrt{dL_k[Tr(R_k^t R_k) + \rho^2(R_k)]/\rho^2(R_k)}},$$

*for $d$ as in lemma 4.3.*

Hence, the estimator is optimal in the sense that the adaptive estimator achieves the best rate of convergence among all the regularized estimators, up to an error of order $\text{pen}(k)$ and $\Sigma(d)/n$. This bound is non asymptotic and the rate of convergence depends on both previous terms.

We point out that the constant $r$ and the weights must be chosen in order to control respectively the penalization term and the constant $\Sigma(d)$. The weights are technically introduced to achieve the so called *Kraft inequality*, $\Sigma(d) < +\infty$ and hence to control the size of the set of parameters $\mathcal{K}_n$.

We also point out that under **SF** $\rho^2(nR_k)$ and $\text{Tr}(R_k^t R_k)/\rho^2(R_k)$ do not depend on $n$.

*Proof.* For any $x_k$ and any $k \in \mathbb{N}$, the mere definition of the estimator $\hat{x}_{\hat{k}}$ implies that

$$\|R_{\hat{k}}(y - T\hat{x}_{\hat{k}})\|^2 + pen(\hat{k}) \leqslant \|R_k(y - Tx_k)\|^2 + pen(k)$$

and

$$\|R_k(y - Tx_k)\|^2 = \|R_k T(x_0 - x_k)\|^2 + 2\langle R_k T(x_0 - x_k), R_k \varepsilon \rangle + \|R_k \varepsilon\|^2$$

Thus, following standard arguments we have

$$\|R_{\hat{k}}T(x_0 - \hat{x}_{\hat{k}})\|^2$$
$$\leqslant \|R_kT(x_0 - x_k)\|^2 - 2 < R_{\hat{k}}T(x_0 - \hat{x}_{\hat{k}}), R_{\hat{k}}\varepsilon >$$
$$+ 2 < R_kT(x_0 - x_k), R_k\varepsilon > -\|R_{\hat{k}}\varepsilon\|^2 + \|R_k\varepsilon\|^2 + \text{pen}(k) - \text{pen}(\hat{k}).$$

Let $0 < \kappa < 1$. Since $2ab \leqslant \kappa a^2 + \frac{1}{\kappa}b^2$, for any $a, b$ we have for any $k$ and $x_k \in X_m$

$$(1-\kappa)\|R_{\hat{k}}T(x_0 - \hat{x}_{\hat{k}})\|^2$$
$$\leqslant (1+\kappa)\|R_kT(x_0 - x_k)\|^2 + \left(2 + \frac{1}{\kappa}\right)\text{pen}(k)$$
$$+ 2\sup_{k \in \mathcal{K}_n}\left\{\left(1 + \frac{1}{\kappa}\right)\|R_k\varepsilon\|^2 - \text{pen}(k)\right\},$$

On the other hand, using that $1 \leqslant R_kT \leqslant C$, we have that for any $x_k \in X_{m_0}$ and any $k \in \mathbb{N}$,

$$(1-\kappa)\|x_0 - \hat{x}_{\hat{k}}\|^2 \leq C(1-\kappa)\|x_0 - x_k\|^2$$
$$+ \left(2 + \frac{1}{\kappa}\right)\text{pen}(k) + 2\,C_1 \sup_{k \in \mathcal{K}_n}\left\{\|R_k\varepsilon\|^2 - \left(1 + \frac{1}{\kappa}\right)^{-1}\text{pen}(k)\right\}.$$

The proof then follows directly from Lemma 4.3 which characterizes the supremum of the empirical process under the linear application as defined by the regularization family. The choice of $\kappa$ will depend on $r$ in the penalty. □

## 3. Applications to standard regularization operators

Penalized estimators are widely used to solve linear inverse problems and can be written in the form (4). Indeed for $k \in \mathcal{K}$, consider a sequence of (diagonal) matrices $A_k \in \mathbb{M}_{m_0 \times m_0}(\mathbb{R})$. Then, define for a chosen sequence $A_k$ the following penalized estimator

$$\hat{x}_k := \arg\min_{x \in X_m}\left[\|\Pi^n_{Y_{m_0}}(y - Tx)\|^2 + \|A_k x\|^2\right]. \tag{7}$$

Note that, for each fixed $A_k$, the expression (7) can be written in the following way:

$$\hat{x}_k = (T^*_{m_0}T_{m_0} + A_k^t A_k)^{-1}T^t_{m_0}y, \tag{8}$$

In practice the second expression is more complicated (the matrix to invert might be big), but it is simpler to deal with in order to show our results concerning the selection of $A_k$. With this notation set the smoothing operator

$$R_k := (T^*_{m_0}T_{m_0} + A_k^t A_k I_{m_0})^{-1}T^t_{m_0}.$$

We point out that choosing the smoothing sequence $A_k$ is the key point since it balances the two terms: if $\|A_k\|$ is large the solution will be smooth but will not, in general, comply to the observations. On the other hand, if $\|A_k\|$ is small, the solution might be too close to the noisy observations to yield a good approximation of $x_0$.

Remark that we can write $R_k = (D^2 + A_k^t A_k I_{m_0})^{-1} D (GG^t)^{-1} G$, as a matrix, with $R_k^t$ its transpose matrix.

### 3.1. Tikhonov regularization

Consider the choice $A_k = \sqrt{\alpha_k} I_{m_0}$, where $\alpha_k$ is a positive decreasing sequence. In this case, the estimator (7) can be written as

$$\hat{x}_k = \arg \min_{x \in X_{m_0}} \left[ \|\Pi_{Y_{m_0}}^n (y - Tx)\|^2 + \alpha_k \|x\|^2 \right]. \quad (9)$$

Where we recognize the usual expression of the Tikhonov estimator. Assume **AN** holds true. Hence Theorem 2.1 can be written in the following way, using the approximation properties of the Tikhonov regularization scheme,

**Proposition 3.1.** *Assume that $\mathcal{K}_n$ is such that $d_{m_0} \geq \sup_{k \in \mathcal{K}_n} \alpha_k^{-1/(2p)}$. Then, by direct calculation of the trace and spectral radius of $R_k$ under* **SF** *and* **SV**, *we have, for $\nu \leq 1$*

$$\frac{Tr(R_k^t R_k)}{\rho^2(R_k)} = O(\alpha_k^{-1/(2p)}).$$

*Also under* **IP** *and* **SC**, *standard approximation results (see for example [11], Chapter 1) the following inequality holds true*

$$\mathbb{E}\|\hat{x}_{\hat{k}} - x_0\|^2 \leq 2\|(I - \Pi_{X_{m_0}})x_0\|^2 + C \inf_{k \in \mathcal{K}_n} [\alpha_k^{2\nu} + 2(\alpha_k^{1+1/2p} n)^{-1}] + \frac{\Sigma(d)}{n}, \quad (10)$$

Hence under **IP** the above inequality yields, for $\nu \leq 1$,

$$\mathbb{E}\|\hat{x}_{\hat{k}} - x_0\|^2 \leq C n^{-\frac{4p\nu}{1+4p\nu+2p}}. \quad (11)$$

Interpreting this rate in the statistical literature reads $s = 2\nu p$: the regularity depends on the ill-posedness of the problem. In the ill-posed literature the error is not related to the underlying dimension so that rates are different. Typically in a Hilbert scale setting, if the true solution $x_0$ belongs to the Hilbert space $H_s$, optimal rates are of order $O(n^{-s/(2s+2p+1)})$, see for example [7]. So, Equation 10 implies that the penalized Tikhonov estimator $\hat{x}_{\hat{k}}$, with $\hat{k}$ selected as in (5), achieves the best rate of convergence over all the possible choices of smoothing sequences $\alpha_k$, $k \in \mathcal{K}_n$. Moreover, if the input data is not too smooth, i.e for $\nu \leqslant 1$, Equation 11 implies that the penalized estimate achieves the minimax rate of convergence for this problem.

To overcome the restrictions induced by the use of Tikhonov estimate, we turn to model selection estimators in the following part.

### 3.2. Regularization by projection

Consider the projection estimator $x_m$ over $X_m$. That is, $x_m$ is chosen in such a way as to minimize $\|\Pi_{Y_m}(y - Tx_m)\|_n$. The estimation error can thus be bounded by

$$\|x_0 - x_m\| \leq \|(I - \Pi_{X_m})x_0\| + \frac{\|\Pi^n_{Y_m}\epsilon\|}{\nu_m}. \tag{12}$$

Since $Y_m$ is a sequence of linear subspaces $\mathbb{E}\|\Pi^n_{Y_m}\epsilon\|^2 = O(\frac{d_m}{n})$. So rates of convergence are of order

$$\|(I - \Pi_{X_m})x_0\| + \frac{d_m^{1/2}}{n^{1/2}\gamma_m}.$$

This rate depends on the ill-posedness of the operator and the approximation properties of $X_m$. Indeed, following [13] and under **SC** we have if $\nu \leq 1/2$,

$$\|(I - \Pi_{X_m})x_0\| \leq \|(I - \Pi^n_{Y_m})T\|^{2\nu} = O(d_m^{-2\nu p}),$$

for a certain $p$. An optimal choice of the dimension (depending on $\nu$) leads to the rate

$$\|x_m - x_0\| = O_{\mathbf{P}}(n^{-\frac{2\nu p}{4\nu p + 2p + 1}}). \tag{13}$$

We aim at using Theorem 2.1 to select the best projection space among a collection of spaces. For this, consider a collection of index sets $m = (j_1, \ldots, j_m)$ such that $m \subset m_0$. For each $m$, define formally $A_m = (a_{ij})$ with

$$\forall i \neq j,\ a_{ij} = 0 \quad \forall i \in m,\ a_{ii} = 0, \quad \forall i \nsubseteq m,\ a_{ii} = \infty.$$

Then we obtain a model selection estimator which can be written as

$$\hat{x}_{\hat{m}} = \arg\min_{(m,x) \in \mathcal{M}_n \times X_m} \|A_{m_0}(y - Tx)\|^2 + \text{pen}(m) \tag{14}$$

where $A_{m_0} = T^+_{m_0}\Pi^n_{Y_{m_0}}$ and $\text{pen}(m)$ is defined as follows. Let $\lambda_{m,j}$ be the eigenvalues of the matrix $\Pi_{X_m}A_{m_0}(\Pi_{X_m}A_{m_0})^t$ and $\{L_m\}_{m \subset m_0}$ a collection of weights such that

$$\Sigma(d) := 2 \sum_{m \subset m_0} \left[\sqrt{\frac{d\sum_{j \in m}\lambda_{m,j}}{\sup_{j \in m}\lambda_{m,j}}} + 1\right] \left[\frac{d}{n\sup_{j \in m}\lambda_{m,j}}\right]^{-1}$$

$$\times e^{-\sqrt{dL_m \frac{\sum_{j \in m}\lambda_{m,j} + \sup_{j \in m}\lambda_{m,j}}{\sup_{j \in m}\lambda_{m,j}}}} < \infty$$

for $d$ as in lemma 4.3. As above, under **SF**, $\Sigma(d)$ does not depend on $n$. Then, for some $\theta > 0$, and $r = 2 + \theta$, set

$$\text{pen}(m) = r(1 + L_m)\sigma^2 \sum_{j \in m} \lambda_{m,j} + \sup_{j \in m} \lambda_{m,j}.$$

**Proposition 3.2.** *Under conditions* **AN, SF, SV, IP** *and* **SC**, *there exists a constant $C$ which depends on $r, T, c_1$ and $c_2$ such that with high probability*

$$\|\hat{x}_{\hat{m}} - x_0\|^2 \leqslant 2\|(I - \Pi_{X_{m_0}})x_0\|^2 + C \inf_{m \in \mathcal{M}_n} [d_m^{-2p\nu} + 2(d_m^{1+2p}n)^{-1}] + \frac{\Sigma(d)}{n}. \quad (15)$$

Hence the penalized model selection estimator is optimal in the sense that it achieves the best rate among the collection of projection estimators. This leads to optimal rates also, as long as there exists a model $m \in \mathcal{M}_n$ with $d_m = n^{1/(4\nu p + 2p + 1)}$, as shown in Equation (13).

We also have the following interpretation for this estimator which offers an important insight. Note that (14) is equivalent to minimizing

$$\begin{aligned}
\hat{x}_{\hat{m}} &= \arg\min_{m \subset m_0} \arg\min_{x_m \in X_m} \{-2 < A_{m_0}y, A_{m_0}Tx_m > + \|A_{m_0}Tx_m\|^2\} + \text{pen}(m) \\
&= \arg\min_{m \subset m_0} \arg\min_{x_m \in X_m} \{-2 < \Pi_{X_m}A_{m_0}y, x_m > + \|x_m\|^2\} + \text{pen}(m).
\end{aligned}$$

Let $\{e_j\}_{j \in m}$ be the canonical base over $X_m$. Define for each $m$,

$$x_{m,j} = < A_{m_0}y, e_j > = < \mathbf{y}, A_{m_0}^t e_j >, \, j = 1, \ldots, m.$$

Thus, $m$ is selected by minimizing

$$-\sum_{j \in m} x_{m,j}^2 + r\sigma^2(1 + L_m)\left[\sum_{j \in m} \lambda_{m,j} + \sup_{j \in m} \lambda_{m,j}\right].$$

We can recognize a hard thresholding scheme, defined in [2], highlighting the equivalence between model selection and penalized M-estimation with an $\ell^0$ penalty, as is also quoted in [17].

## 4. Appendix

In this section we give some technical lemmas. The next lemma characterizes the supremum of an empirical process by the norm of an orthogonal projection.

**Lemma 4.1.**
$$\sup_{y \in Y_m, \, \|y\|_n = 1} | < \varepsilon, y >_n | = \|\Pi_{Y_m}^n \varepsilon\|_n \quad (16)$$

*Proof.* Using the definition of an orthogonal projector, we have

$$\left\|\varepsilon - \frac{1}{\|\Pi_{Y_m}^n \varepsilon\|_n} \Pi_{Y_m}^n \varepsilon\right\|_n^2 = \min_{\{y \in Y_m, \, \|y\|_n = 1\}} \|\varepsilon - y\|_n^2.$$

As a consequence we can write:

$$\|\varepsilon\|_n^2 - 2 < \varepsilon, \frac{1}{\|\Pi_{Y_m}^n \varepsilon\|_n} \Pi_{Y_m}^n \varepsilon >_n + \frac{1}{\|\Pi_{Y_m}^n \varepsilon\|_n^2} \|\Pi_{Y_m}^n \varepsilon\|_n^2$$
$$= \min_{\{y \in Y_m, \, \|y\|_n = 1\}} \|\varepsilon\|_n^2 - 2 < \varepsilon, y >_n + 1$$

$$2 < \varepsilon - \Pi^n_{Y_m}\varepsilon, \frac{1}{\|\Pi^n_{Y_m}\varepsilon\|_n}\Pi^n_{Y_m}\varepsilon >_n + 2 < \Pi^n_{Y_m}\varepsilon, \frac{1}{\|\Pi^n_{Y_m}\varepsilon\|_n}\Pi^n_{Y_m}\varepsilon >_n$$

$$= 2 \sup_{\{y \in Y_m,\, \|y\|_n=1\}} | <\varepsilon, y>_n |$$

$$\|\Pi^n_{Y_m}\varepsilon\|_n = \sup_{\{y \in Y_m,\, \|y\|_n=1\}} | <\varepsilon, y>_n |,$$

which ends the proof. □

The next result is a deviation inequality based on a functional exponential inequality (Theorem 7.4) due to [4] 2003.

**Lemma 4.2.** *Set $\eta(A) = \sup_{\|u\|=1} \sum_{i=1}^n \varepsilon_i (A^t u)_i$ for $A : \mathbb{R}^n \to \mathbb{R}^k$. Let*

$$v = \mathbb{E}\sum_{i=1}^n \sup_{\|u\|=1} \frac{(A^t u)_i^2}{\rho(A^t A)}\left(\frac{\varepsilon_i}{\sigma}\right)^2 + 2\mathbb{E}\eta(A)/(\sigma\rho^{1/2}(A^t A)).$$

*Then,*

$$P\left(\frac{\eta(A)}{\sigma\rho^{1/2}(A^t A)} > \mathbb{E}\frac{\eta(A)}{\sigma\rho^{1/2}(A^t A)} + \sqrt{2vx} + x\right) \leq e^{-x}.$$

*Proof.* Since the application $u \to A^t u$ is continuous, we have $\eta(A) = \sup_{u \in S} \sum_{i=1}^n \varepsilon_i (A^t u)_i$ for $S$ some countable subset of the unit ball. On the other hand,

$$\sup_{\|u\|=1} [A^t u]_i \leq \sup_{\|u\|=1} \|A^t u\| \leq \rho(A).$$

Thus $\sup_{\|u\|\leq 1} |(A^t u)_i / \rho^{1/2}(A^t A)| \leq 1$. Also, following the proof of Corollary 5.1 in [1]

$$\sup_{\|u\|=1} \frac{(A^t u)_i^2}{\rho(A^t A)}$$
$$\leq \sup_{\|u\|=1} \frac{(\sum_{j=1}^m u_j (A^t e_j)_i)^2}{\rho(A^t A)}$$
$$\leq \sup_{\|u\|=1} \frac{(\sum_{j=1}^m (A^t e_j)_i^2) \sum_{j=1}^m u_j^2}{\rho(A^t A)}$$
$$:= z_i.$$

Set $Z = Z(\varepsilon_1, \ldots, \varepsilon_n) = \eta(A)/(\sigma\rho^{1/2}(A^t A))$. Let $\mathbb{E}_j$ stand for the conditional expectation given $\varepsilon_i$ for $i \neq j$. Hence, in the proof of Theorem 7.4 in [4] we may bound

$$|Z - \mathbb{E}_j Z|^q \leq \frac{|\varepsilon_j|^q}{\sigma^q} \sup_{\|u\|=1} \frac{(A^t u)_j^2}{\rho(A^t A)} \sup_{\|u\|=1} \max_i \left(\frac{(A^t u)_i^2}{\rho(A^t A)}\right)^{q-2} \leq (|\varepsilon_j|/\sigma)^q z_j.$$

Thus, $\mathbb{E}|Z - \mathbb{E}_j Z|^q \leq z_j q!/2$. Finally, note that

$$\sum_{j=1}^n z_j = \frac{Tr(A^t A)}{\rho(A^t A)}.$$

Thus, the proof follows from Theorem 7.4 in [4]. □

As a corollary, we have the following lemma

**Lemma 4.3.** • *There exists a positive constant d that depends on r such that the following inequality holds*

$$P(\eta^2(A) \geq \sigma^2[Tr(A^tA) + \rho(A^tA)]r/2(1+L) + \sigma^2 u) \quad (17)$$
$$\leq \exp\{-\sqrt{d(1/\rho(A^tA)u + r/2L[Tr(A^tA)/\rho(A^tA) + 1])}\}$$

- *Set $k_1 = d/(\rho(A^tA)\sigma^2)$ and $k_2 = dr/2L[Tr(A^tA)/\rho(A^tA)+1]$. Then, there exists a constant $C_q$, which depends only on $q$, such that,*

$$\mathbb{E}[\eta^2(A) - \sigma^2[Tr(A^tA) + \rho(A^tA)]r/2(1+L)]_+^q \quad (18)$$
$$\leq C_q k_1^{-q}[k_2^{q-1/2} + k_2^{q-1}]e^{-\sqrt{k_2}}$$

*holds.*

*Proof.* As a first step we will bound $v$. Since $\mathbb{E}Z \leq \mathbb{E}^{1/2}Z^2$, we have

$$\begin{aligned} v &\leq \mathbb{E}\sum_{i=1}^n z_i\left(\frac{\varepsilon_i}{\sigma}\right)^2 + 2\sqrt{\mathbb{E}\sum_{i=1}^n z_i\left(\frac{\varepsilon_i}{\sigma}\right)^2} \\ &\leq (1+\nu)\mathbb{E}\sum_{i=1}^n z_i\left(\frac{\varepsilon_i}{\sigma}\right)^2 + Tr(A^tA)/\rho(A^tA). \end{aligned}$$

Moreover, following, [1] p. 480, for all $q \geqslant 2$, the following version of Rosenthal's inequality holds:

$$\mathbb{E}^{q/2}\sum_{i=1}^n z_i\left(\frac{\varepsilon_i}{\sigma}\right)^2 \leq 2^{q/2}Tr(A^tA)/\rho(A^tA)\mathbb{E}\frac{|\varepsilon_1|^q}{\sigma^q}.$$

Hence, we have

$$v \leq (1+\nu)Tr(A^tA)/\rho(A^tA) + \frac{1}{\nu}$$

and

$$v^2 \leq 2\left[2^2(1+\nu)^2 Tr(A^tA)/\rho(A^tA)\mathbb{E}\frac{|\varepsilon^4|}{\sigma^4} + \left(\frac{1}{\nu}\right)^2\right].$$

Set $0 < \alpha < 1$. Choose $\delta$ and $\beta$ such that if

$$2^2 4!\delta^2(1+1/\alpha)(1-\nu)^2 < c_1,$$
$$2\delta^2(1+1/\alpha)(\frac{1}{\nu})^2 < c_2$$

and $c = \max((1+\beta)\max(c_1,c_2),(1+\beta)(1+\alpha))$, then $r/2 > c$. Let $u > 0$ and without loosing generality, assume $\sigma = 1$. Thus,

$$P(\eta^2(A) \geq (Tr(A^tA) + \rho(A^tA))r/2(1+L) + u)$$
$$\leq P(\eta^2(A) \geq (Tr(A^tA)(1+\alpha) + (1+1/\alpha)\delta^2\nu^2\rho(A^tA))(1+\beta)$$
$$+[r/2 - c](Tr(A^tA) + \delta^2v^2\rho(A^tA)) + r/2L(Tr(A^tA) + \rho(A^tA)) + u)$$
$$\leq P(\eta^2(A) \geq (Tr(A^tA)(1+\alpha) + (1+1/\alpha)\delta^2\nu^2\rho(A^tA))(1+\beta)$$
$$+r/2L(Tr(A^tA) + \rho(A^tA)) + u)$$

Set
$$x' = \left(1 + \frac{1}{\beta}\right)^{-1}\left[\frac{r}{2L}\left(\frac{Tr(A^tA)}{\rho(A^tA)} + 1\right) + \frac{u}{\rho(A^tA)}\right].$$

The last term is equal to
$$P\left(\frac{\eta^2(A)}{\rho(A^tA)} \geq \left(\frac{Tr(A^tA)}{\rho(A^tA)}(1+\alpha)\right.\right.$$
$$\left.\left. + (1+1/\alpha)v^2\delta^2\right)(1+\beta) + r/2L\frac{Tr(A^tA)}{\rho(A^tA)} + 1\right) + u\right)$$
$$= P\left(\frac{\eta^2(A)}{\rho(A^tA)} \geq \left(\frac{Tr(A^tA)}{\rho(A^tA)}(1+\alpha) + (1+1/\alpha)v^2\delta^2\right)(1+\beta) + (1+1/\beta)x'\right)$$

Finally, we may bound
$$\leq P\left(\frac{\eta^2(A)}{\rho(A^tA)} \geq \left(\mathbb{E}\frac{\eta(A)}{\rho^{1/2}(A^tA)} + \delta v\right)^2(1+\beta) + (1+1/\beta)x'\right)$$
$$\leq P\left(\frac{\eta^2(A)}{\rho(A^tA)} \geq \left(\mathbb{E}\frac{\eta(A)}{\rho^{1/2}(A^tA)} + \delta v + \sqrt{x'}\right)^2\right)$$
$$= P\left(\frac{\eta(A)}{\rho^{1/2}(A^tA)} \geq \mathbb{E}\frac{\eta(A)}{\rho^{1/2}(A^tA)} + \delta v + (1+2/\delta)x''\right)$$
$$\leq P\left(\frac{\eta(A)}{\rho^{1/2}(A^tA)} \geq \mathbb{E}\frac{\eta(A)}{\rho^{1/2}(A^tA)} + \sqrt{2vx''} + x''\right) \leq e^{-x''}$$
$$= e^{-\sqrt{g(A)}},$$

where we have used repeatedly that for any constant $c > 0, ca^2 + 1/cb^2 \geq 2ab$ and set

$$g(A) = ((1+1/\beta)^{-1}(1+2/\delta)^2)(r/2L[Tr(A^tA)/\rho(A^tA) + 1] + u/\rho(A^tA)).$$

Set also $d = [(1+1/\beta)^{-1}(1+2/\delta)^2]^{-1}$ and $b(A) = Tr(A^tA)/\rho(A^tA)$. Thus we have shown the first part of the lemma.

Moreover, using the above inequality,

$$\mathbb{E}[\eta^2(A) - \sigma^2(Tr(A^tA) + \rho(A^tA))r/2(1+L)]^q_+$$
$$\leq \int_0^\infty \sigma^{2q}qu^{q-1}e^{-\sqrt{dr/2L[b(A)+1]+du/(\rho(A^tA))}}du.$$

Consider the change of variable $w = du/(\rho(A^tA)) + dr/2L[b(A) + 1]$, so that

$$\mathbb{E}[\eta^2(A) - \sigma^2(Tr(A^tA) + \rho(A^tA))r/2(1+L)]_+^q$$
$$\leq \left(\frac{\sigma^2\rho(A^tA)}{d}\right)^q \int_{dr/2L[b(A)+1]}^{\infty} (w - dr/2L[b(A) + 1])^{q-1} e^{-\sqrt{w}} dw.$$

The last expression is in turn bounded by

$$\left(\frac{\sigma^2\rho(A^tA)}{d}\right)^q \int_{dr/2L[b(A)+1]}^{\infty} e^{-\sqrt{w}}[w^{q-1} + (dr/2L[b(A) + 1])^{q-1}] dw$$
$$\leq C_q k_1^{-q}[k_2^{q-1/2} + k_2^{q-1}]e^{-\sqrt{k_2}},$$

ending the proof.

$\square$